\documentclass{aims}
\usepackage{amsmath}
  \usepackage{paralist}
  \usepackage{graphics} 
  \usepackage{epsfig} 
\usepackage{graphicx}  \usepackage{epstopdf}
\input cyracc.def

\usepackage{filecontents}

 \usepackage[colorlinks=true]{hyperref}
\hypersetup{urlcolor=blue, citecolor=red}
\RequirePackage[hyperpageref]{backref} 
    \renewcommand*{\backref}[1]{}  
    \renewcommand*{\backrefalt}[4]{
       \ifcase #1 
          No cited.
       \or
          Cited on p. #2.
       \else
          Cited on pp. #2.
       \fi} 

\def\currentvolume{X}
 \def\currentissue{X}
  \def\currentyear{200X}
   \def\currentmonth{XX}
    \def\ppages{X--XX}
     \def\DOI{10.3934/xx.xx.xx.xx}

\newtheorem{theorem}{Theorem}[section]

\newtheorem{lemma}[theorem]{Lemma}

\theoremstyle{definition}
\newtheorem{definition}[theorem]{Definition}
\newtheorem{remark}{Remark}

\newcommand{\ep}{\varepsilon}
\newcommand{\eps}[1]{{#1}_{\varepsilon}}

    \def\BbbE{\mathbb E}
    \def\BbbN{\mathbb N}

\def\bH{\mathbf H}
\def\bG{\mathbf G}
\def\bP{{\bf P}}
\def\bE{{\bf E}}
\def\cF{{\mathcal F}}
\def\cB{{\mathcal B}}
\def\cT{{\mathcal T}}

\def\cS{{\mathcal S}}
\def\cK{{\mathcal K}}

\def\frakP{\frak P}
\def\vecPi{\overrightarrow{\Pi}}

\def\one{{\mathbb I}}
\newcommand{\ones}{1\!\!\one}

\def\vecalpha{{\overrightarrow\alpha}}
        
        \def\vecpi{{\overrightarrow{\Pi}}}
        \def\vecPi{\overrightarrow{\Pi}}
        \def\uPi{\underline{\Pi}}
        
\def\ux{\underline{x}}

\def\uA{\underline{A}}
\def\uB{\underline{B}}
\def\uD{\underline{D}}
\def\uX{\underline{X}}
\def\uf{\underline{f}}
\def\udelta{\underline{\delta}}
\def\ualpha{\underline{\alpha}}
\def\plus{+} 
\newcommand*{\esssup}{\mathop{\mathrm{ess\,sup}}\displaylimits}

\title[Detection of distributed disorders] 
      {A Model of Distributed Disorders Detection}

\author[Krzysztof J.  Szajowski]{}


 \email{Krzysztof.Szajowski@pwr.edu.pl}

\thanks{$^*$ Corresponding author: Krzysztof J. Szajowski}

\begin{document}
\maketitle

\centerline{\scshape Krzysztof J. Szajowski$^*$}
\medskip
{\footnotesize
 \centerline{Wrocław University of Science and Technology}
   \centerline{Faculty of Pure and Applied Mathematics}
   \centerline{Wybrzeże Wyspiańskiego 27, 50-370 Wrocław, Poland}
} 

\medskip


\bigskip

 \centerline{(Communicated by the associate editor Arnd R\"osch)}

\begin{abstract}
The paper deals with disorders detection in the multivariate stochastic process. We consider the  multidimensional Poisson process or the multivariate renewal process. This class of processes can be used as a description of the distributed detection system. The multivariate renewal process can be seen as the sequence of random vectors, where parts of its coordinates are holding times, others are the size of jumps and the index of stream, at which the new event appears. It is assumed that at each stream two kinds of changes are possible: in the holding time or in the size of jumps distribution. The various specific mutual relations between the change points are possible. The aim of the research is to derive the detectors which realize the optimal value of the specified criterion. The change point moment estimates have been obtained in some cases. The difficulties have appeared for the dependent streams with unspecified order of change points. The presented results suggest further research on the construction of detectors in the general model.
\end{abstract}

\section{\label{intro} Introduction}
The subject of discussion is the model that describes the phenomenon of piecewise deterministic  signals. Time-intervals between jumps are random variables and jumps' size is also random. The modeled object, for a random period, is in a homogeneous state. In a random time the signal changes its nature and time between jumps, although further random, but have a different distribution, or the size of jumps changes its distribution. After the change is a time homogeneous process and at random time another change. In the present case, there may be one change in the distribution of time between notifications and one change in the distribution of jumps' size. The aim is to locate the two changes in real time. The other change in the behavior of the process does not continue to alter. Signals of this nature appear in technical issues, medicine, and finance. 

\subsection{A preliminary consideration} 
In this paper the construction of the mathematical model of the phenomenon described above requires the determination of a probabilistic space in which all random quantities, random variables and stochastic processes, are defined. Let $(\Omega, \cF, \bP(\cdot))$ be fixed probability space. The consideration will focus on the renewal-reward model of change points detection. A~\emph{renewal--reward process} is a jump process with a general holding time distributions and  general distribution of jumps (see Br\'{e}maud~\cite{bre81:point}, Jacobsen~\cite{jac06:point}). Let us denote $\{W_i\}_{i=1}^\infty$ a sequence of \emph{iid} \emph{rv} (\emph{rewards}) with the finite expected value. The random variable $Y_t = \sum_{i=1}^{X_t}W_i$, where $X_t$ is the renewal process, is called a renewal-reward process. On the turn, a \emph{renewal process} is a pure jump process with a general distribution of the holding times. To describe this type of processes, let us consider a sequence of positive \emph{iid} \emph{rv} $\{ S_i\}_{i=1}^\infty$ with the finite and positive expectation.  The renewal process $X_t = \sum^{\infty}_{n=1} \mathbb{I}_{\{J_n \leq t\}}=\sup \left\{\, n: J_n \leq t\, \right\}$, where  $J_n = \sum_{i=1}^n S_i$  for each $n > 0$. 
The mentioned changes of distributions are assumed to appear on the sequences  $\{ S_i\}_{i=1}^\infty$ and $\{W_i\}_{i=1}^\infty$ at some random moments $\theta_1$ and $\theta_2$, respectively. The general detectors of disorders take values at $\Re^{+}$. In this approach the class of decision function will be restricted to the event (the jump) moments. This means that the moments of the decision can be identified with the index jumping moments in the process. The underlined models can be represented as the sequences of random vectors $\{(S_n,W_n)\}_{n=1}^\infty$. The \emph{r.v.} $W_n$ and $S_n$ do not need to be necessarily independent in contrast to the classical theory of the renewal--reward processes. The aim of the research is to formulate the rigorous model of the problem and to investigate them. 

\subsection{A motivation}
The motivation for the project is the wide literature related to the optimal stopping problem and the scarce results for the multiple stopping settings. The related questions for the change point models are also stimulating. There is satisfactory literature describing the state of art in the disorder problem or the change-point problems in the off-line and on--line setting. Let us mention the monograph by Brodsky and Darkhovsky \cite{brodar93:nonparametr} and the story by Shiryaev \cite{shi06:filtring}. 

The disorder problem has a long history (see Shiryaev~\cite{shi06:filtring}, Proceedings of AMS-IMS-SIAM Summer Research Conference on {''Change point problems''}  E.~Carlstein et al. (ed)~\cite{carmulsie94:hazard}). It is known from the early papers by Page~\cite{pag54:inspection,pag55:change}, Girshick and Rubin~\cite{girrub52:Bayes}. However, it was A.N.~Kolmogorow who had formulated the realistic and mathematically precise model of rapid detection  of a change point or disorder. His student, A.N.~Shiryaev~\cite{shi61:quickest,shi63:disorder} (see the history described in \cite{shi06:filtring}) published the first important results for the sequential problems of disorder detection in the Poisson distribution. This direction of research is the subject of intensive work of contemporary projects. Let us mention papers by Gal{}’chuk and Rozowski~\cite{galroz71MR0297028:Poisson}, Davis~\cite{dav74:Poisson}, Peskir and Shiryaev~\cite{pesshi02:Poisson} (see also Gapeev~\cite{Gap05MR2115049:Poisson}, Bayraktar, Dayanik and Karatzas \cite{bayerhday06:Poisson}.
The manifold experiment to formulate the adequate model of disorder for more complex processes which appear in the observed phenomena of nature and economy, showed extreme difficulties (see Fuh~\cite{fuh04:hidden}, Ivanoff and Merzbach~\cite{ivamer10:change-set}, Szajowski~\cite{Sza11:multivariate}). 

The work continues to carry out research described in the papers published by Sarnowski and the author~\cite{SarSza08:disorder}, \cite{SarSza11:transition} on the change point problem for the undefined, the Markovian type processes before and after the disorder. Related results can be found in the papers by Bojdecki i Hosza~\cite{bojhos84:problem}, Szajowski~\cite{sza92:detection}, 
Yoshida~\cite{yos83:complicated}, Yakir~\cite{yak94:finite} and  Moustakides~\cite{mou98:abrupt}. 
The most important guidelines for considering the results of this work are contained in the authors' works~\cite{Sza11:random} on multivariate disorders detection. 

The key technique used in the work is based on a multiple optimal stopping a Markov process. The fundamental knowledge on the optimal stopping of random processes can be found in the monographies by Chow, Robbins and Siegmund~\cite{chosie71:exp}, \cite{shi08:MR2374974} or Peskir and Shiryaev~\cite{pesshi06:optimalfree}. The first work devoted to the multiple stopping of the discrete time sequences processes was published by Haggstroma~\cite{hag67:twostop} and for the Markov processes Nikolaev~\cite{nik81:obob} (see also Eidukjavicjus~\cite{eid79:dvumia}, Nikolaev\cite{nik98:multstopping}). The extension of the model to semi-Markov processes has been made by Stadje~\cite{sta85:multiple}. Recent results for continuous time processes have been presented by Kobylanski, Quenez and  Rouy-Mironescu~\cite{kobquerou10:double}) and the various applications of such an approach are described in papers by Feng and Xiao \cite{fenxia00:continuous}, \cite{fenxi00:optimal}, Karpowicz and Szajowski~\cite{karsza07:pric,karsza12:anglers}. For the risk process such extension can be found in papers by Karpowicz and Szajowski \cite{karsza07:risk}. However, in the approach applied here there is no need to refer to the results for the semi-Markov processes. The analytical difficulties for such problems open various questions concerning the construction of algorithms, approximate solutions and Monte Carlo methods.

\subsection{Variation of the basic problems}
A general formulation starts from a random sequences $\{ S_i\}_{i=1}^\infty$ and $\{W_i\}_{i=1}^\infty$ having segments that are the homogeneous Markov sequences. Each segment has its own transition probability law and the length of the segment is unknown and random. It means that there are random moments $\theta_1$ and $\theta_2$, at which the source of observations is changed. The transition probabilities of each process at the segment are chosen from a given set of distributions. An \emph{a priori} distributions of the disorder moments are given.  
The problem is to construct the detection algorithm of the disorders. The algorithm should detect the change points with the given precision and maximal probability. 

To be more precise, let us see examples.  It is easy to ascertain that in the renewal--reward process the distribution of holding time can be first changed in the moment $\theta_1$, and next, in the  moment $\theta_2\geq \theta_1$ the distribution of the random variables $W_n$ changes. One can formulate three such problems:
\begin{itemize}
\item it is known that $\theta_1<\theta_2$, \emph{ie} holding time has the distribution change moment earlier than the rewards distribution;
\item the order of change points is different: the first the reward sequences are disordered;
\item there is no information which sequences will change distribution as the first.
\end{itemize}
In Bayesian model the last case needs additional a priori information about the chance that the disorder of rewards will appear before the holding time disorder. 

If the process is a model of signals from the \emph{distributed sources} then the disorders in different sources are combined to construct the final decision about the reason of non-homogeneous behavior of the process. Based on the idea of a simple game the model of the fusion center is proposed. The strategy of detection at each segment (source) of the process is defined as the equilibrium in a non-cooperative game between the selfish sensors (see \cite{Sza11:multivariate}).
   
For the single process having structure similar to those of the generalized compound Poisson processes the temporal disorder is possible. It is a natural problem, mentioned by researchers (see \cite{BayPoo07:quickest}, \cite{yos83:complicated} and others). When the model of the process is equivalent to the sequence of a vector of random variables, it is possible that each coordinate changes their distribution in some moments which could be different at each coordinate. 

Our environment is described by the states of nature and they are viable over time. There are no objective boundaries in the state space defining safe or unsafe world. These borders are marked by a subjective knowledge. After the appointment of the task of maintaining a safe environment the aim is to observe and monitor the states in time to know their relation to the boundaries. Methodology presented in the paper allows to prepare the environment description and the detectors of the borders of the areas. For the illustration of a description it was selected variant environment in which states form the renewal--reward process. Areas of the state space, acceptable or not, describe the observed distributions of the process. Acceptable states of the system have a distribution, and unacceptable otherwise. Due to the multi-dimensional character of the process the specified areas are determined by the boundaries of the components of the state vector. The model is analysed and the boundaries are determine for the areas of warning states. To determine the boundaries of the detector indicates that to have the correct detection it is necessary a priori knowledge of the possible scenarios achieve a state of emergency. Efficient detection of threats is only possible if the nature implements provided by our scenario. Strategy to little knowledge about the possible scenario is far less effective and it is much more cumbersome to implement.

\section{Random switching between Markov processes}
Let us formulate the general detection for at most two change points problem or switching detection problem. et us consider an observable sequence of random variables $(X_n)_{ n \in \BbbN}$  with the  homogeneous structure on the time intervals: $n\in [0,\theta_1-1]$, $[\theta_1,\theta_2-1]$ and $[\theta_2,\infty)$. The parameters $\theta_1$, $\theta_2$ are pairs of random variables with values in $\BbbN$ having distributions:
\begin{eqnarray}
\label{rozkladyTeta}
\bP(\theta_1 = j) &=& \one_{\{j=0\}}(j)\pi+\one_{\{j>0\}}(j)\bar{\pi} p_1^{j-1}q_1,\\
\label{rokladWarTeta2}
\bP(\theta_2 = k \mid \theta_1=j) &=&\one_{\{k=j\}}(k)\rho+\one_{\{k>j\}}(k)\bar{\rho} p_2^{k-j-1}q_2
\end{eqnarray}
where $j=0,1,2,...$, $k=j,j+1,j+2,...$, $\bar\pi=1-\pi$, $\bar\rho=1-\rho$. 

In segments the distribution depends additionally on the parameter $\epsilon_i$ with values from the finite set $\cK=\{1,2,\ldots,K\}$. It is assumed that $\bP(\epsilon_i = s)=r^{i}_s$, $s\in\cK$ and $i\in\{0,1,2\}$.

In the sequel if  there is only one possible process in the segment then the second index will be abandoned. 

Additionally, there are Markov processes $(X_n^{is}, \mathcal{G}_n^{is}, \bP_x^{is})$, $i=0,1,2$, $s\in\cK$, where $\sigma$-fields $\mathcal{G}_n^{is}$ are the smallest $\sigma$-fields for which
 ${(X^{is}_n)}_{n=0}^\infty$ are adapted, respectively. The process $(X_n)_{ n \in \BbbN}$ is  connected with the random variables $\theta_1$, $\theta_2$, $\epsilon_i$ and the Markov processes $\{X_n^{is}\}_{n=0}^\infty$  as follows:
\begin{eqnarray}
\label{procesyX}
X_n &=& X^{0s_0}_n \one_{\{\theta_1>n,\epsilon_0=s_0\}} + X^{1s_1}_{n-\theta_1+1}\one_{\{X^{1}_{0}=x^{0}_{\theta_1-1},\theta_1 \leq n < \theta_2,\epsilon_1=s_1\}}\\
&&\mbox{} + X^{2s_2}_{n-\theta_2+1}\one_{\{X^{2}_{0}=x^{1}_{\theta_2-\theta_1},\theta_2 \leq n,\epsilon_2=s_2\}}
\nonumber .
\end{eqnarray}
The observable sequence of \emph{rv} is defined on the space $(\Omega,\mathcal{F}, \bP)$ with values in Borel subset $(\BbbE, \mathcal{B})$, $\BbbE\subset \mathbf{R}$ with $\sigma$-additive measure $\mu$. The measures ${\bP^i_x}(\cdot)$ on $\cF$, $i=0,1,2$, have the following representation: 
\begin{eqnarray*}
\bP^{is}_x(\omega:X_1^{is}\in B)&=&\int_B f_x^{is}(y)\mu(dy)=\int_B \mu_x^{is}(dy)=\mu^{is}_x(B),
\end{eqnarray*}
for any $B\in\cB$, where $f_x^{is}(\cdot)$ are different and $f_x^{is_i}(y)/f_x^{((i+1)\text{mod} 3)s_{(i+1)\text{mod} 3}}(y) < \infty$ for $i=0,1,2$, $s_\cdot\in\cK$ and all $x,y \in \BbbE$. 
\subsection{Finite dimensional distribution of process}
For any $D_n=\{\omega:X_i\in B_i,\; i=1,\ldots,n\}$, where $B_i\in \cB$, and any $x\in\BbbE$ define
\begin{eqnarray*}
\bP_x(D_n)&=&\int_{\times_{i=1}^n B_i}S_n(x,\vec{y}_n) \mu(d\vec{y}_n)=\int_{\times_{i=1}^nB_i} \mu_x(d\vec{y}_n)=\mu_x(\times_{i=1}^nB_i). 
\end{eqnarray*} 
Let $\cS$ be the set of all stopping times with respect to $(\cF_n)$, $n=0,1,\ldots$ and
$\cT=\{(\tau,\sigma): \tau\leq \sigma,\ \tau,\sigma\in\cS\}$. 

\subsection{Criteria of change--point detection}
The aim of the DM is to indicate the moments of switching with given precision $d_1,d_2$ (Problem $\mbox{D}_{d_1d_2}$). We want to determine a pair of stopping times $(\tau^*,\sigma^*)\in\cT$
such that for every $x\in\BbbE$
\begin{eqnarray}
\label{equ3}
&&\bP_x(-d_{1,l}\leq\tau^*-\theta_1\leq d_{1,r},d_{2,l}\leq\sigma^*-\theta_2\leq d_{2,r})\\
\nonumber&&\mbox{\;} =
\sup_{\stackrel{(\tau,\sigma)\in\cT}{0\leq\tau\leq\sigma<\infty}} 
\bP_x(d_{1,l}\leq\tau-\theta_1\leq d_{1,r},d_{2,l}\leq\sigma-\theta_2\leq d_{2,r}).
\end{eqnarray}

The problem with the fixed transition distributions at each segment has been formulated in \cite{sza96:twodis} and has been extended to the case $0\leq \theta_1\leq \theta_2$ in \cite{Sza11:random}. 
The investigated models here assume that between $\theta_1$ and $\theta_2$, the distribution is chosen from a given set (for simplification having two elements only). The distribution is predetermined in two models and chosen randomly in one model.    

{Let us introduce the following notation:}
\begin{eqnarray*}
L_{st}^{u_1}(\underline{x}_{k,n}) &=&\displaystyle{ \prod\limits_{r=k+1}^{n-s-t}} f_{x_{r-1}}^{0}(x_{r})\displaystyle{ \prod\limits_{r=n-s-t+1}^{n-t}}f_{x_{r-1}}^{1u_1}(x_{r})\displaystyle{ \prod\limits_{r=n-t+1}^{n}}f_{x_{r-1}}^{1}(x_{r}), \nonumber\\
\underline{A}_{k,n} &=&  \times_{i=k}^{n} A_{i} = A_k \times A_{k+1} \times \ldots \times A_n,\; \mbox{ $A_i \in \cB$, $u_1\in\cK$} \nonumber
\end{eqnarray*}
where the convention $\prod_{i=j_1}^{j_2}x_i = 1$ for $j_1 > j_2$ is used.

Let $B_i \in {\cB}$, $m\leq i \leq n$ and let us assume that $X_0=x$ and denote $\uD_{m,n}=\{\omega:X_i(\omega)\in A_i, m\leq i\leq n\}$. For $D_i=\{\omega:X_i\in A_i\}\in\cF_i$, $m\leq i\leq n$ we have by properties of the density function $S_n(\ux_{1,n})$ with respect to the measure $\mu(\cdot)$ 
\begin{eqnarray*}
\int_{A_{s,t}}L_{m,n}^{u_1}(\ux_{s-1,t})\mu(d\ux_{s,t})=\frakP_{m,n}^{u_1}(X_{s-1},\uD_{s,t}),
\end{eqnarray*}
where $m+n\leq t-s+1$, $u_1\in\cK$.
Let us now define functions $S_\cdot(\cdot)$ and $H_\cdot(\cdot,\cdot,\cdot,\cdot)$ and the sequence of functions $S_n:\times_{i=1}^n\BbbE\rightarrow\Re$ as follows: $S_0(x_0)=1$ and for $n\geq 1$:
\setlength\arraycolsep{0pt}
\small
\begin{eqnarray}\label{disXuncond}
S_n(\ux_n)
&=&f^{\epsilon_1,\theta_1\leq\theta_2\leq n}_x(\ux_{1,n})+f^{\epsilon_1,\theta_1\leq n<\theta_2}_x(\ux_{1,n})\\
\nonumber&&\mbox{}+f^{\epsilon_1,\theta_1=\theta_2>n}_x(\ux_{1,n})+f^{\epsilon_1,n<\theta_1<\theta_2}_x(\ux_{1,n}).
\end{eqnarray}
Aditionally, we have
\begin{eqnarray}\label{disHuncond}
\bH(\cdot,\cdot,\cdot,\cdot)&=&f(x_{n+1}|x_n).
\end{eqnarray}
For further calculation it is important to have $n$-dimensional distribution for various configurations of disorders. 
\begin{eqnarray}
\label{theta1lesstheta2lessn} f^{\epsilon_1,\theta_1\leq\theta_2\leq n}_x(\ux_{1,n}) &=&\bar{\pi}\rho\sum_{u\in\cK}r^1_u\{\sum_{j=1}^{n}p_1^{j-1}q_1 L_{0,n-j+1}^{u}(\ux_{0,n})\\
\nonumber &\hspace{-18em}+&\hspace{-9em}\bar{\pi}\bar{\rho} \sum_{j=1}^{n-1}\sum_{k=j+1}^{n}\{p_1^{j-1}q_1p_2^{k-j-1}q_2L_{k-j,n-k+1}^{u}(\ux_{0,n})\}+\pi\rho L_{0,n}^{u}(\ux_{0,n})\}\\
\label{theta1less_n_lesstheta2} f^{\epsilon_1,\theta_1\leq n<\theta_2}_x(\ux_{1,n}) &=&\bar{\rho}\sum_{u\in\cK}r^1_u[\bar{\pi}\sum_{j=1}^n\{p_1^{j-1}q_1 p_2^{n-j}L_{n-j+1,0}^{u}(\ux_{0,n})\}\\ \nonumber &&\mbox{}+\pi\sum_{j=1}^n\{p_2^{j-1}q_2L_{j-1,n-j+1}^{u}(\ux_{0,n})\}]\\
\label{theta1=theta2gen} f^{\epsilon_1,\theta_1=\theta_2>n}_x(\ux_{1,n})&=& \rho\bar{\pi}p_1^n\sum_{u\in\cK}r^1_uL_{0,0}^{u}(\ux_{0,n})\\
\label{nlesstheta1lesstheta2} f^{\epsilon_1,n<\theta_1<\theta_2}_x(\ux_{1,n})&=& \bar{\rho}\bar{\pi}p_1^n\sum_{u\in\cK}r^1_uL_{0,0}^{u}(\ux_{0,n}).
\end{eqnarray}
Denote $\left\langle\ \underline{u}\:,\underline{v}\right\rangle = \sum_{i=1}^d u_i v_i$  and $\ones=(1,1,\ldots,1)\in\Re^K$.
We have (cf. \cite{SarSza08:disorder}, \cite{dubmaz01:quickest})
\begin{lemma}\label{recSn}
For $n>0$ the function $S_n(\ux_{1,n})$ follows recursion
\begin{align}
\label{recurSn}S_{n+1}(\ux_{1,n+1})&=\bH(x_n,x_{n+1},\vecPi_n,\Upsilon_n)S_n(\ux_{1,n})\\
\intertext{where for $\vec{\pi}=(\alpha,\beta,\gamma)$, and $\upsilon=(\upsilon^1,\upsilon^2,\ldots,\upsilon_K)$\hfil}
\label{Hfunction}\bH(x,y,\vec{\pi},\upsilon)&= (1-\alpha) p_1 f^{0}_x(y)+ [q_2\alpha+p_2\beta + q_1\gamma]f^{2}_x(y)\\
\nonumber&\mbox{}+  \langle\ p_2(\alpha -\beta)+ q_1(\upsilon- \alpha-\gamma)\:,\underline{f^{1}_x}(y)\rangle .
\end{align}
\end{lemma}
Here $\vecPi_n=(\underline{\Pi_n^1},\underline{\Pi_n^2},\underline{\Pi_n^{12}})$ and $\Upsilon_n=(\Upsilon_n^1,\Upsilon_n^2,\ldots,\Upsilon_n^K)$.
One can formulate: $\bH(x,y,\vec{\pi},\upsilon)=\left\langle\ \ones \:, \underline{\bH}\right\rangle =\sum_{u\in\cK}\underline{\bH}^u(x,y,\vec{\pi}^u,\upsilon^u)$.
Let us assume $0\leq \theta_1\leq \theta_2$ and suppose that $B_i \in {\cB}$, $1\leq i \leq n+1$ and  $X_0=x$ and denote $D_n=\{\omega:X_i(\omega)\in B_i, 1\leq i\leq n\}$.
For $A_i=\{\omega:X_i\in B_i\}\in\cF_i$, $1\leq i\leq n+1$ we have by properties of the density function $S_n(\ux_{1,n})$ with respect to the measure $\mu(\cdot)$ 
\begin{eqnarray*}
\int_{D_{n+1}}d\bP_x&=&\int_{D_n}\one_{A_{n+1}}d\bP_x.
\end{eqnarray*}

Now we split the conditional probability of $\{X_{n+1}\in A_{n+1},\epsilon_1=u\}$ into the following parts
\begin{eqnarray}
\nonumber\bP_{x}(X_{n+1}\in A_{n+1},\epsilon_1=u  \mid \cF_n ) \\ \label{nProb1}&\hspace{-12em}=&\hspace{-5em}\bP_{x}(n<\theta_1<\theta_2,X_{n+1} \in A_{n+1},\epsilon_1=u \mid \cF_n) \\
&\hspace{-12em}+&\hspace{-5em}\bP_{x}(\theta_1 \leq n < \theta_2,X_{n+1} \in A_{n+1},\epsilon_1=u \mid \cF_n) \label{nProb2}\\
&\hspace{-12em}+&\hspace{-5em} \bP_{x}(n < \theta_1 =\theta_2,X_{n+1} \in A_{n+1},\epsilon_1=u \mid \cF_n) \label{nProb3}\\
&\hspace{-12em}+&\hspace{-5em} \bP_{x}(\theta_1\leq\theta_2 \leq n,X_{n+1} \in A_{n+1},\epsilon_1=u \mid \cF_n)\label{nProb4}
\end{eqnarray}
\begin{description}
\item[In \eqref{nProb1}] we have:
\begin{eqnarray*}
&\hspace{-8em}\mbox{\;}&\hspace{-3em}\int_{D_n}\bP_{x}(\theta_2>\theta_1 > n,X_{n+1} \in A_{n+1},\epsilon_1=u \mid \cF_n)d\bP_x\\
&\hspace{-12em}=&\hspace{-1em} \int_{\times_{i=1}^{n}B_i}\hspace{-2em} (f^{\epsilon_1 = u,n<\theta_1<\theta_2}_x(\ux_{1,n})\int_{B_{n+1}}(p_1f^0_{x_n}(x_{n+1})+q_1f^{1,u}_{x_n}(x_{n+1}))\mu(dx_{n+1}))\mu(d\ux_{1,n})\\
&\hspace{-12em}=&\hspace{-1em}\int_{D_n}\bP_{x}(\theta_2>\theta_1 > n,\epsilon_1 = u \mid \cF_n) [ \bP_{X_n}^0(A_{n+1})p_1+q_1\bP_{X_n}^{1,u}(A_{n+1}) ] d\bP_x.
\end{eqnarray*}
\end{description}
\begin{description}
\item[In \eqref{nProb2}] we get by similar arguments as for (\ref{nProb1})
\footnotesize
\begin{eqnarray*}
&\hspace{-20em}\mbox{\;}&\hspace{-10em} \bP_{x}(\theta_1 \leq n < \theta_2,X_{n+1} \in A_{n+1},\epsilon_1 = u \mid \cF_n) \nonumber\\
 &=& \left( \bP_{x}(\theta_1 \leq n,\epsilon_1 = u \mid \cF_n) - \bP_{x}(\theta_2 \leq n,\epsilon_1 = u \mid \cF_n) \right)\nonumber\\
 && \times \! [ q_2 \bP_{X_n}^2(A_{n+1}) + p_2\bP_{X_n}^{1,u}(A_{n+1})] \nonumber
\end{eqnarray*}

\item[In \eqref{nProb4}] this part has the form:
\[
\bP_{x}(\theta_2 \leq n, X_{n+1} \in A_{n+1},\epsilon_1 = u\mid \cF_n) = \bP_{x}(\theta_2 \leq n,\epsilon_1 = u \mid \cF_n)\bP_{X_n}^2(A_{n+1})
\]
\item[In \eqref{nProb3}] the conditional probability is equal to
\begin{eqnarray*}
&\hspace{-20em}\mbox{\;}&\hspace{-10em}\bP_{x}(\theta_1=\theta_2>n , X_{n+1} \in A_{n+1},\epsilon_1 = u \mid \cF_n) \nonumber\\
 &=& \bP_{x}(\theta_1 =\theta_2> n,\epsilon_1 = u \mid \cF_n)[ q_1 \bP_{X_n}^2(A_{n+1}) + p_1\bP_{X_n}^{0}(A_{n+1})] \nonumber
\end{eqnarray*}
\end{description}
These formulae lead to
\[
f(X_{n+1}|\uX_{1,n})=\bH(X_n,X_{n+1},\Pi_n^1,\Pi_n^2,\Pi_n^{12},\Upsilon_n).
\]
\subsection{$n$-dimensional vs. $n+d$-dimensional distributions}
It is not too difficult to get the recursive formula for $\bH(\cdot)$.
\begin{align}
\label{Hufunction}
\bH(x,y,\vec{\pi},\upsilon)&=\sum_{u\in\cK}\underline{\bH}^u(x,y,\vec{\pi}^u,\upsilon^u)\\
\underline{\bH}^u(x,y,\vec{\pi}^u,\upsilon^u)&= (\upsilon^u-\alpha^u) p_1 f^{0}_x(y)+  [p_2(\alpha^u -\beta^u)+ q_1(\upsilon^u- \alpha^u-\gamma^u)]f^{1,u}_x(y) \\
\nonumber &  + [q_2\alpha^u+p_2\beta^u + q_1\gamma^u]f^{2}_x(y).
\end{align}
\begin{lemma}\label{recSn}
For $n>0$ the density function $S_n(\ux_{1,n})$ follows recursion
\begin{align}
\label{recurSn}S_{n+d}(\ux_{1,n+d})&=\bG_d(\ux_{n,n+d},\vecPi_n,\Upsilon_n)S_n(\ux_{0,n})\\
\intertext{where\hfil}\vspace{-1.5ex}
\label{Hfunction} \bG_d(\ux_{n,n+d},\vecPi_n,\Upsilon_n)&=f(\ux_{n+1,n+d}\mid \ux_{0,n})
\end{align}
\end{lemma}
Now we split the conditional probability of $\uA_{n+1,n+d}$ into the following parts
\begin{eqnarray*}\label{n+dProb1}
\bP_{x}(\uX_{n+1,n+d}\in \uB_{n+1,n+d}  \mid \cF_n )\\
&\hspace{-14em}=&\hspace{-7em} \sum_{u\in\cK}[\bP_{x}(n<\theta_1<\theta_2,\uX_{n+1,n+d} \in \uB_{n+1,n+d},\epsilon_1 = u \mid \cF_n) \\
&\hspace{-12em}+&\hspace{-5em} \bP_{x}(\theta_1 \leq n < \theta_2,\uX_{n+1,n+d} \in \uB_{n+1,n+d},\epsilon_1 = u \mid \cF_n) \label{n+dProb2}\\
&\hspace{-12em}+&\hspace{-5em} \bP_{x}(n < \theta_1 =\theta_2,\uX_{n+1,n+d} \in \uB_{n+1,n+d},\epsilon_1 = u \mid \cF_n) \label{n+dProb3}\\
&\hspace{-12em}+&\hspace{-5em} \bP_{x}(\theta_1\leq\theta_2 \leq n,\uX_{n+1,n+d}\in \uB_{n+1,n+d},\epsilon_1 = u  \mid \cF_n)]\label{n+dProb4}
\end{eqnarray*}

\begin{eqnarray*}
&\hspace{-12em}\mbox{\;}&\hspace{-6em}\int_{D_n}\bP_{x}(\theta_2>\theta_1 > n,\uX_{n+1,n+d} \in \uB_{n+1,n+d},\epsilon_1 = u \mid \cF_n)d\bP_x\\
&=&\int_{D_n}\bP_{x}(\theta_2>\theta_1 > n,\epsilon_1 =u \mid \cF_n) [ \sum_{j=1}^dp_1^{j-1}q_1\frakP_{j,d-j+1}^u(X_n,\uA_{n+1,n+d})\\
&&+p_1^d\frakP_{0,0}^u(X_n,A_{n+1,n+d}) ] d\bP_x.
\end{eqnarray*}
Let us calculate: $\bP_{x}(\theta_1 \leq n < \theta_2,\uX_{n+1,n+d} \in \uB_{n+1,n+d},\epsilon_1 = u \mid \cF_n)$. This probability can be calculated as follows:
\begin{eqnarray*}
\mbox{} &=& \sum_{j=1}^d\bP_{x}(\theta_1 \leq n < \theta_2, \theta_2 = n+j, X_{n+1,n+d} \in \uB_{n+1,n+d},\epsilon_1 = u \mid \cF_n) \nonumber \\
 &&+ \bP_{x}(\theta_1 \leq n < \theta_2, \theta_2 \neq n+d, \uX_{n+1,n+d} \in \uB_{n+1,n+d},\epsilon_1 = u \mid \cF_n) \nonumber \\
 &=& \left( \bP_{x}(\theta_1 \leq n ,\epsilon_1 = u\mid \cF_n) - \bP_{x}(\theta_2 \leq n,\epsilon_1 = u \mid \cF_n) \right)\nonumber\\
 && \times \! [\sum_{j=1}^d p_2^j \frakP_{j,d-j+1}^u(X_n,\uA_{n+1,n+d}) + p_2^d\bP_{X_n}^{1,u}(A_{n+1})] \nonumber
\end{eqnarray*}

\section{On some a posteriori processes}
For \eqref{equ3} the following \emph{a posteriori} processes are crucial (cf. \cite{yos83:complicated}, \cite{sza92:detection}).
\begin{eqnarray}
\label{pi1ux}
\Pi^{i,u}_n&=&\bP_x(\theta_i\leq n,\epsilon_1 = u|\cF_n),\\
\label{pi1x}
\Pi^i_n&=&\left\langle\ \ones \:, \underline{\Pi^{i}_n}\right\rangle= \sum_{u\in\cK}\Pi^{i,u}_n=\bP_x(\theta_i\leq n|\cF_n),\\
\label{pi12ux}
\Pi^{12,u}_{n}&=&P_x(\theta_1=\theta_2>n,\epsilon_1 = u|\cF_{nn}),\\
\label{pi12nmx}
\Pi^{12}_{n}&=&\left\langle\ \ones \:, \underline{\Pi^{12}_n}\right\rangle =\bP_x(\theta_1=\theta_2>n|\cF_n),\\
\label{pinmux}
\Pi_{mn}^{u}&=&\bP_x(\theta_1=m,\theta_2>n,\epsilon_1 = u|\cF_{mn}),\\
\label{pinmx}
\Pi_{mn}&=&\left\langle\ \ones \:, \underline{\Pi^{i}_{mn}}\right\rangle =\bP_x(\theta_1=m,\theta_2>n|\cF_{mn}),\\
\label{upsnmx}
\Upsilon_n^u&=&\bP_x(\epsilon_1 = u|\cF_{n}),
\end{eqnarray}
for $m,n=1,2,\ldots$, $m<n$, $i=1,2$. Also $\cF_{n}=\cF_{nn}$.
\subsection{Recursive form of posteriors}
A posteriori processes: ${\Pi}^{1,u}_{n}$, ${\Pi}^{2,u}_{n}$, ${\Pi}^{12,u}_{n}$, $\uPi_{m\;{n}}^u$ can be calculated based on the following formula:
\begin{eqnarray}
\label{pi1xu}
\overline{\underline{\Pi}}^{1,u}_{n+1}&=&\overline{\uPi}^{1,u}_{n}\frac{p_1f_{x_n}^{0}(x_{n+1})}{\bH(x,y,\uPi^{1}_{n},\uPi^{2}_{n},\uPi^{12}_{n},\Upsilon_n)}\\
\label{pi2xu}
\uPi^{2,u}_{n+1}&=&\frac{(q_2\uPi^{1,u}_{n} + p_2\uPi^{2,u}_{n} + q_1\uPi^{12,u}_{n}) f^{2}_x(y)}{\bH(x,y,\uPi^{1}_{n},\uPi^{2}_{n},\uPi^{12}_{n},\Upsilon_n)},\\
\label{pi12nmx}
\uPi^{12,u}_{n+1}&=&\frac{p_1\uPi^{12,u}_{n} f^0_x(y)}{\bH(x,y,\uPi^{1}_{n},\uPi^{2}_{n},\uPi^{12}_{n},\Upsilon_n)},\\
\label{pinmx}
\uPi_{m\;{n+1}}^u&=&\frac{p_2\uPi_{m\;{n}}^u f^{1,u}_x(y)}{\bH(x,y,\uPi^{1}_{n},\uPi^{2}_{n},\uPi^{12}_{n},\Upsilon_n)},
\end{eqnarray}
for $m,n=1,2,\ldots$, $m<n$, $i=1,2$.

\subsection{Markov processes based on observations} 
Let us definie $\Upsilon_{n}^u$ in the recursive way.
\begin{eqnarray}
\label{Upsilon1n}
\Upsilon_n^u&=&\Pi^{1,u}_n+\overline{\Pi}^{1,u}_n,\\
\label{Upsilon1n1}
\Upsilon_{n+1}^u&=&\frac{f_{x_n}^{0}(x_{n+1})[\Upsilon_{n}^u-q_1\uPi^{12,u}_{n}]+f^{1,u}_x(y)p_2[\uPi^{1,u}_{n}-\uPi^{2,u}_{n}]}{\bH(x,y,\uPi^{1}_{n},\uPi^{2}_{n},\uPi^{12}_{n},\Upsilon_n)}\\
\nonumber&&+\frac{f_{x_n}^{2}(x_{n+1})[q_2(\uPi^{1,u}_{n}+\uPi^{2,u}_{n})+\uPi^{2,u}_{n}+q_1\uPi^{12,u}_{n}]}{\bH(x,y,\uPi^{1}_{n},\uPi^{2}_{n},\uPi^{12}_{n},\Upsilon_n)}
\end{eqnarray}
For recursive representation of (\ref{pi1xu})--(\ref{Upsilon1n1}) we need the following functions:
\begin{eqnarray*}
\label{pi1}
\uPi^{1,u}(x,y,\alpha,\beta,\gamma,\upsilon)&=&\upsilon^u-[p_1(\upsilon^u-\alpha^u) f^{0}_x(y)]{\bH^{-1}(x,y,\alpha,\beta,\gamma,\upsilon)}\\
\label{pi2}
\uPi^{2,u}(x,y,\alpha,\beta,\gamma,\upsilon)&=&[(q_2\alpha^u + p_2\beta^u + q_1\gamma^u) f^{2}_x(y)]{\bH^{-1}(x,y,\alpha,\beta,\gamma,\upsilon)}\\
\label{pi12}
\uPi^{12,u}(x,y,\alpha,\beta,\gamma,\upsilon)&=&p_1\gamma^u f^0_x(y){\bH^{-1}(x,y,\alpha,\beta,\gamma,\upsilon)}\\
\label{pinm}
\uPi^u(x,y,\alpha,\beta,\gamma,\delta,\upsilon)&=&p_2\delta^u f^{1,u}_x(y){\bH^{-1}(x,y,\alpha,\beta,\gamma,\upsilon)}.
\end{eqnarray*}

The function $\underline{\bH}^u$ and $\bH$ will be used in representations of the posteriors.
\begin{eqnarray*}
\underline{\bH}^u(x,y,\alpha,\beta,\gamma,\upsilon)&=& (\upsilon^u-\alpha^u) p_1 f^{0}_x(y)\\
&&\mbox{}+ [p_2(\alpha^u -\beta^u) + q_1(\upsilon^u- \alpha^u-\gamma^u)]f^{1,u}_x(y)\\
&&\mbox{}  + [q_2\alpha^u+p_2\beta^u + q_1\gamma^u]f^{2}_x(y).
\end{eqnarray*} 
There is $\bH(x,y,\alpha,\beta,\gamma,\upsilon)=\left\langle\ \ones\:,\underline{\bH}\right\rangle$.
{In the sequel we adopt the following denotations: $\vecalpha=(\alpha,\beta,\gamma)$ and $\vecpi_n^u=(\uPi^{1,u}_n,\uPi^{2,u}_n,\uPi^{12,u}_n)$.}

\section{The disorder problem vs. the stopping problem} 
The basic formulae used in the transformation of the disorder problems to the stopping problems are given in the following
\begin{lemma}
\label{reqform}
For each $x\in\BbbE$ and each Borel function $u: \Re\rightarrow\Re$ the following formulae  for $m,n=1,2,\ldots$, $m<n$, $i=1,2$, hold:
\begin{eqnarray}
\label{eqpi1}
 \Pi^{i,u}_{n+1}&=&\uPi^{i,u}(X_n,X_{n+1},\uPi^1_n,\uPi^2_n,\uPi^{12}_n,\Upsilon_n)\\
 \label{eqpi12}
 \Pi^{12,u}_{n+1}&=&\uPi^{12,u}(X_n,X_{n+1},\uPi^1_n,\uPi^2_n,\uPi^{12}_n,\Upsilon_n)\\
 \label{eqpi}
 \Pi_{m\,n+1}^u&=&\uPi(X_n,X_{n+1},\uPi^1_n,\uPi^2_n,\uPi^{12}_n,\uPi_{m\,n},\Upsilon_n)
\end{eqnarray}
with the boundary condition $\Pi^{1,u}_0=\pi r_u$, $\Pi^{2,u}_0(x)=\pi r_u\rho$,
$\uPi_{m\,m}^u=(1-\rho)\frac{q_1 f^{1,u}_{X_{m-1}}(X_m)}{p_1 f^0_{X_{m-1}}(X_m)}(\Upsilon_n^u-\uPi^{1,u}_m)$.

\end{lemma}
\begin{lemma}
\label{multidistform}
For the problem ~\ref{equ3} the following formulae are valied: 
\begin{enumerate}
\item $\bP_x(\theta_2=\theta_1>n+1,\epsilon=u|\cF_n)=p_1\uPi^{12,u}_n$;
\item $\bP_x(\theta_2>\theta_1>n+1,\epsilon=u|\cF_n)=p_1(\Upsilon_n^u-\uPi_n^{1,u}-\uPi_n^{12,u})$;
\item $\bP_x(\theta_1\leq n+1,\epsilon=u|\cF_n)=\bP(\theta_1\leq n+1<\theta_2,\epsilon=u|\cF_n)+\bP(\theta_2\leq n+1,\epsilon=u|\cF_n)$;
\item $\bP(\theta_1\leq n+1<\theta_2,\epsilon=u|\cF_n)=q_1(1-\uPi^{1,u}_n-\uPi^{12,u}_n)+p_2(\uPi^{1,u}_n-\uPi^{2,u}_n)$;
\item $\bP_x(\theta_\leq n+1,\epsilon=u|\cF_n)=q_2\uPi^{1,u}_n+p_2\uPi_n^{2,u}+q_1\uPi^{12,u}_n$.
\end{enumerate}
\end{lemma}
\begin{lemma}
    \label{lematTechniczny}
For each $x\in\mathbf{E}$ and each Borel function $u:\mathbf{R}\longrightarrow \mathbf{R}$ the following equations are fulfilled.
\begin{eqnarray}
\label{wartOczPi1}
\bE_{x}\left(u(X_{n+1})(\Upsilon_n^u - \uPi_{n+1}^{1,u})\mid \mathcal{F}_n\right) &=& (1-\uPi^{1,u}_n-\uPi^{12,u}_n)p_1\\
\nonumber&&\times\int_{\BbbE}u(y)f^0_{X_n}(y)\mu_{X_n}(dy),\\
\nonumber\bE_{x}\left(u(X_{n+1})(\uPi_{n+1}^{1,u} - \uPi_{n+1}^{2,u})\mid \cF_n\right)&=&
 \left[q_1(1-\uPi^{1,u}_n-\uPi^{12,u}_n)+p_2(\uPi^{1,u}_n-\uPi^{2,u}_n) \right]\\
\label{wartOczPi1Pi2}&&\times    \int_{\BbbE}u(x)f^{1,u}_{X_n}(y)\mu_{X_n}(dy),\\
\label{wartOczPi2}
\bE_{x}\left(u(X_{n+1})\uPi_{n+1}^{2,u})\mid \cF_n\right) &=&
\left[ q_2\uPi^{1,u}_n+p_2\uPi_n^{2,u}+q_1\uPi^{12,u}_n \right]\\
\nonumber&&\times\int_{\BbbE}\!u(y)f^2_{X_n}(y)\mu_{X_n}(dy),\\
\label{wartOczPi12}
\bE_{x}\left(u(X_{n+1})\uPi^{12,u}_{n+1})\mid \cF_n\right) &=&
\left[p_1\uPi^{12,u}_n \right]\!\int_{\BbbE}\!u(y)f^0_{X_n}(y)\mu_{X_n}(dy)\\
\label{exp1}
\bE_x(u(X_{n+1})|\cF_n) &=&
\int_{\BbbE} u(y)\bH(X_n,y,\vecpi_n(x))\mu_{X_n}(dy)
\end{eqnarray}
\end{lemma}
\section{An equivalent issue -- the double optimal stopping problem}\label{dosp}
A {\it compound stopping variable} is a pair $(\tau,\sigma)$ of stopping times such that $\tau\leq\sigma$ a.e.. Denote $\cT_m=\{(\tau,\sigma)\in\cT: \tau\geq m\}$, $\cT_{mn}=\{(\tau,\sigma)\in\cT: \tau=m, \sigma\geq n\}$ and $\cS_m=\{\tau\in\cS:\tau\geq m\}$ ($\cF_{mn}=\cF_n$, $m,n\in\BbbN$, $m\leq n$). 

We define two-parameter stochastic sequence 
$\xi(x)=\{\xi_{mn}\}_{m,n\in\BbbN,\  m<n,\ x\in\BbbE}$, where $\xi_{mn}=\bP_{x}(\theta_1=m, \theta_2=n|\cF_{mn})$.

\subsection{The optimal compound stopping variable}
For every $x\in\BbbE$, $m,n\in\BbbN$, $m<n$ let us define the optimal stopping problem of $\xi(x)$ on 
$\cT^\plus_{mn}=\{(\tau,\sigma)\in\cT_{mn}:\tau<\sigma\}$.  
A compound stopping variable $(\tau^*,\sigma^*)$ is said to be optimal in $\cT^\plus_m$ if $\bE_{x}\xi_{\tau^*\sigma^*} = \sup_{(\tau,\sigma)\in\cT^\plus_m}\bE_{x}\xi_{\tau\sigma}$ \\
(or in $\cT^\plus_{mn}$ if $\bE_{x}\xi_{\tau^*\sigma^*} =
\sup_{(\tau,\sigma)\in\cT^\plus_{mn}}\bE_{x}\xi_{\tau\sigma}$). 
Let us define
\begin{equation}
\label{pr1}
\eta_{mn} =
\esssup_{(\tau,\sigma)\in\cT^\plus_{mn}}\bE_x(\xi_{\tau\sigma}|\cF_{mn}).
\end{equation}
If we put $\xi_{m\infty}=0$, then 
\[
\eta_{mn}=\esssup_{(\tau,\sigma)\in\cT^\plus_{mn}}\bP_{x}(\theta_1=\tau,\theta_2=\sigma|\cF_{mn}).
\]
From the theory of optimal stopping for the double indexed processes (cf. \cite{hag67:twostop}, \cite{nik81:obob}) the sequence $\eta_{mn}$ satisfies 
\[
\eta_{mn}=\max\{\xi_{mn},\bE(\eta_{mn+1}|\cF_{mn})\}.
\]
If $\sigma^*_m=\inf\{n>m: \eta_{mn}=\xi_{mn}\}$, then
$(m,\sigma^*_n)$ is optimal in $\cT^\plus_{mn}$ and $\eta_{mn}=\bE_x(\xi_{m\sigma^*_n}|\cF_{mn})$ a.e.. 
Define
$\hat{\eta}_{mn}=\max\{\xi_{mn},\bE(\eta_{m\;n+1}|\cF_{mn})\}$  for $n\geq m$ 
if $\hat{\sigma}^*_m=\inf\{n\geq m: \hat{\eta}_{mn}=\xi_{mn}\}$, then $(m,\hat{\sigma}^*_m)$ is optimal in $\cT_{mn}$ and $\hat{\eta}_{mm}=\bE_x(\xi_{m\sigma^*_m}|\cF_{mm})$ a.e.. 
\begin{lemma}
\label{lem1}
Stopping time $\sigma^*_m$ is optimal for every stopping problem (\ref{pr1}).
\end{lemma}
What is left is to consider the optimal stopping problem for $(\eta_{mn})_{m=0,n=m}^{\infty,\infty}$ on $(\cT_{mn})_{m=0,n=m}^{\infty,\infty}$. For further consideration denote $\eta_{m} = \bE_{x}(\eta_{mm+1}|\cF_m)$.
The first stop payoff will be determined. Let us define:
\begin{equation}
\label{FirstStopPayment}
V_m=\esssup_{\tau\in\cS_m}\bE_x(\eta_\tau|\cF_m).
\end{equation}
Then $V_m=\max\{\eta_m,\bE_x(V_{m+1}|\cF_m)\} \mbox{ a.e.}$ and
we define $\tau^*_n=\inf\{k\geq n: V_k=\eta_k\}$.
\begin{lemma}[\ref{lem2}]
\label{lem2}
The strategy $\tau^*_0$ is the optimal strategy of the first stop.
\end{lemma}
\subsection{Solution of the equivalent double stopping problem}
For this presentation the case $d_1=d_2=0$ is considered. Let us construct multidimensional  Markov chains such that $\xi_{mn}$ and $\eta_m$ will be the functions of their states. By considering the \emph{a posteriori} processes we get $\xi_{00}=\pi\rho$ and for $m<n$
\begin{equation}
\label{ximm}
\xi_{m\,n}^x \stackrel{L. \ref{multidistform}}{=} \bP_x(\theta_1=m,\theta_2=n|\cF_{m\,n})
=\left\{\begin{array}{ll} \frac{q_2}{p_2}\langle\ \uPi_{m\,n}(x)\:,\frac{f^2_{X_{n-1}}(X_n)}{\uf^1_{X_{n-1}}(X_n)}\rangle&\text{ for $m<n$}\\
\rho\frac{q_1}{p_1}\frac{f^2_{X_{m-1}}(X_m)}{f^0_{X_{m-1}}(X_m)}(1-\Pi_m^1)&
\text{for $n=m$.}
\end{array}
\right.
\end{equation}
 
The vector $(X_n,X_{n+1},\vecpi_n,\uPi_{m\,n},\Upsilon_n)$ for $n=m+1, m+2,\ldots$ is a function of $(X_{n-1},X_{n},\vecpi_{n-1},\uPi_{m\,{n-1}},\Upsilon_n)$ and
$X_{n+1}$. Besides the conditional distribution of $X_{n+1}$
given $\cF_n$ depends on $X_n$, $\uPi^1_n(x)$,  $\uPi^2_n(x)$ and $ \Upsilon_n$ only. 

These facts imply that  
$\{(X_n,X_{n+1},\vecpi_n,\uPi_{m\,n},\Upsilon_n )\}_{n=m+1}^\infty$
form a homogeneous Markov process. This allows us to reduce the basic problem
(\ref{pr1}) for each $m$ to the optimal stopping problem of the Markov process
$Z_m(x)=\{(X_{n-1},X_n,\vecpi_n,\uPi_{m\,n},\Upsilon_n),\ m,n\in\BbbN,\ 
 m<n,\ x\in\BbbE\}$ with the reward function $h(t,u,\vecalpha,\delta,\upsilon)= \frac{q_2}{p_2}\langle\udelta,\frac{f^2_t(u)}{\uf^1_t(u)}\rangle$. 
\begin{lemma}
 \label{lem4}
 A solution of the optimal stopping problem (\ref{pr1}) for $m=1,2,\ldots$ has a form
\begin{equation*}
\label{optsstop}
\sigma^*_m=\inf\{n>m:
\langle\frac{\uPi_{m\; n}}{\Pi_{m\; n}}, \frac{f^2_{X_{n-1}}(X_n)}{\uf^1_{X_{n-1}}(X_n)}\rangle\geq R^*(X_n,\uPi_{m\; n})\}
\end{equation*}
where $R^*(t,\udelta)=p_2\int_{\BbbE}\langle r^*(t,s,\udelta),\uf^1_t(s)\rangle\mu_t(ds)$ and the function $r^*(t,u,\udelta)$ satisfies the equation $r^*(t,u,\udelta)=\max\{\langle\frac{\udelta}{\delta},\frac{f^2_t(u)}{\uf^1_t(u)}\rangle,p_2\int_{\BbbE} \langle r^*(u,s,\udelta),\uf^1_u(s)\rangle\mu_u(ds)\}$. The value of the problem is equal
\begin{equation*}
\label{valpr1}
\eta_m = \bE_x(\eta_{m\,{m+1}}|\cF_m)      =\frac{q_1}{p_1}\langle\frac{\uf^1_{X_{m-1}}(X_m)}{f^0_{X_{m-1}}(X_m)}, \overline{\uPi}^1_m\rangle R_\rho^\star(X_{m-1},X_m,\uPi_{m\; m}),
\end{equation*}
where $R_\rho^\star(t,u,\udelta)=\max\{\rho\langle\frac{\udelta}{\delta},\frac{f^2_{t}(u)}{\uf^1_{t}(u)}\rangle,\frac{q_2}{p_2}(1-\rho)R^\star(t,\udelta)\}$.
\end{lemma}
Based on the results of Lemma~\ref{lem4} and properties of the \emph{a posteriori} process $\Pi_{nm}$ we have the optimal second moment 
\[
\hat\sigma_0^\star=\left\{\begin{array}{ll}
0&\mbox{if $\pi\rho\geq q_1 (1-\pi)\int_\BbbE f_x^1(u)R^\star_\rho(x,u,\udelta)\mu_x(du)$,}\\
\sigma_0^\star&\mbox{ otherwise.}
\end{array}
\right.
\]
By lemmas \ref{lem4} and \ref{reqform} (formula \eqref{eqpi}) the optimal stopping problem
\eqref{FirstStopPayment} has been transformed to the optimal stopping
problem for the homogeneous Markov process $W=\{(X_{m-1},X_m,\vecpi_m,\Pi^{12}_m,\Upsilon_m),\ m\in\BbbN,\ x\in\BbbE\}$
with the reward function
\begin{equation*}
\label{FstepPayoff}
f(t,u,\vecalpha,\upsilon)=\frac{q_1}{p_1}\langle\frac{\uf^1_{t}(u)}{f^0_{t}(u)},\bar{\alpha}\rangle R^\star_\rho(t,u,\udelta).
\end{equation*}
\begin{theorem}
 \label{lem5}
 A solution of the optimal stopping problem \eqref{FirstStopPayment} for $n=1,2,\ldots$ has a form
\begin{equation}
\label{optfstop}
\tau^*_n=\inf\{k\geq n: (X_{k-1},X_k,\vecpi_k,\upsilon_k)\in B^*\}
\end{equation}
where $B^*=\{(t,u,\vecalpha,\upsilon):\langle \bar{\ualpha}, \frac{f^2_t(u)}{\uf^1_t(u)} \rangle R^\star_\rho(t,u,\udelta)\geq
p_1\int_{\BbbE} v^*(u,s,\vecalpha,\upsilon)f^0_u(s)\mu_u(ds)\}$.

The function
$v^*(t,u,\vecalpha,\upsilon)=\lim_{n\rightarrow\infty} v_n(t,u,\vecalpha,\upsilon)$, where
$v_0(t,u,\vecalpha,\upsilon)=R^\star_\rho(t,u,\udelta)$, 
\begin{eqnarray*}
v_{n+1}(t,u,\vecalpha,\upsilon)&=&\max\{\langle\bar{\ualpha},\frac{f^2_t(u)}{f^1_t(u)}\rangle R^\star_\rho(t,u,\udelta),\\
&&p_1\int_{\BbbE}v_n(u,s,\vecalpha,\upsilon)\langle \bar{\ualpha},\uf^1_u(s)\rangle\mu_u(ds)\}.
\end{eqnarray*}
$v^*(t,u\vecalpha,\upsilon)$ satisfies the equation

The value of the problem $V_n=v^*(X_{n-1},X_n,\vecPi_n,\Upsilon_n)$.
\end{theorem}

\bibliography{\jobname}
\bibliographystyle{abbrvnat}

\medskip
Received 31st of March 2019; revised 31st of January 2020.
\medskip

\end{document}